\documentclass[12pt]{article}

\usepackage{tikz, enumitem, caption, subcaption, float}
\usepackage{indentfirst, hyperref}
\usepackage{amsmath,amsfonts,amssymb,amsthm,epsfig,epstopdf,titling,url,array, mathdots}
\usepackage[margin=1in]{geometry}
\usepackage{mathtools}

\usepackage{nccmath}
\newtagform{roman}[]()

\usetikzlibrary{snakes}
\usetikzlibrary{decorations.markings}
\tikzset{unlabelled/.style={black, draw=black, circle, fill, scale=0.5}}
\tikzset{labelled/.style={black, draw=black, circle, scale=0.7, font=\Large}}

\tikzset{->-/.style={decoration={
			markings,
			mark=at position #1 with {\arrow{>}}},postaction={decorate}}}
\tikzset{dashedarc/.style={decoration={
			markings,
			mark=at position #1 with {\arrow{>}}},postaction={decorate}}}
		
\usepackage{algorithm, algorithmic}

\begin{document}

\newtheorem{theorem}{\hspace{5mm}Theorem}[section]
\newtheorem{corollary}[theorem]{\hspace{5mm}Corollary}
\newtheorem{lemma}[theorem]{\hspace{5mm}Lemma}
\newtheorem{proposition}[theorem]{\hspace{5mm}Proposition}
\newtheorem{definition}[theorem]{Definition\hspace{5mm}}

\newcommand{\lemmaautorefname}{Lemma}
\newcommand{\corollaryautorefname}{Corollary}
\newcommand{\propositionautorefname}{Proposition}

\newcommand{\pf}{{\bf Proof: }}

\newenvironment{AMS}{}{}
\newenvironment{keywords}{}{}

\title{Obstructions for acyclic local tournament orientation completions}

\author{Kevin Hsu\ and\ Jing Huang 
     \thanks{Department of Mathematics and Statistics,
      University of Victoria, Victoria, B.C., Canada; huangj@uvic.ca}}
\date{}

\maketitle

\begin{abstract}
The orientation completion problem for a fixed class of oriented graphs asks
whether a given partially oriented graph can be completed to an oriented graph in
the class. Orientation completion problems have been studied recently for several
classes of oriented graphs, yielding both polynomial time solutions and 
NP-completeness results.
Local tournaments are a well-structured class of oriented graphs that generalize 
tournaments and their underlying graphs are intimately related to proper 
circular-arc graphs. Proper interval graphs are precisely those which can
be oriented as acyclic local tournaments. It has been proved that the orientation 
completion problems for local tournaments and acyclic local tournaments are both
polynomial time solvable. In this paper we identify the obstructions for acyclic
local tournament orientation completions. These are in a sense minimal partially 
oriented graphs that cannot be completed to acyclic local tournaments. 
Our description of the obstructions imply that they can be recognized 
in polynomial time. In a companion paper we will determine all obstructions for 
local tournament orientation completions.
\end{abstract}

\section{Introduction}

Many graph properties can be reformulated in terms of the existence of certain 
orientations. A celebrated theorem of Robbins \cite{robbinsAMM46} states that 
a graph is 2-edge-connected if and only if it has a strong orientation.
Gallai, Roy and Vitaver \cite{gallai1968,royRFIRO1,vitaverDANS147}
proved that a graph is $k$-colourable if and only if it has an orientation 
in which the longest directed path has at most $k$ vertices.

An oriented graph $D = (V,A)$ is called {\em transitive} if for any three vertices
$u, v, w$, $(u,v) \in A$ and $(v,w) \in A$ imply $(u,w) \in A$. The underlying 
graphs of transitive oriented graphs are the {\em comparability graphs} 
\cite{gallai}. 
This is equivalent to say that a graph is a comparability graph if and only if  
it has a transitive orientation. An oriented graph $D = (V,A)$ is called 
{\em quasi-transitive} if for any three vertices $u, v, w$, $(u,v) \in A$ and 
$(v,w) \in A$ imply $(u,w) \in A$ or $(w,u) \in A$, cf. \cite{bh}. Transitive 
oriented graphs are all quasi-transitive and the converse is not true. 
Interestingly, the underlying graphs of quasi-transitive oriented graphs are also 
the comparability graphs, cf. \cite{gh}.

A graph $G = (V,E)$ is a {\em proper circular-arc graph} if there is a family of 
circular-arcs $I_v, v \in V$ on a circle such that no circular-arc 
contains another and for all $u, v \in V$, $uv \in E$ if and only if $I_u$ and $I_v$
intersect. It is easy to see that every proper circular-arc graph has an orientation
$D = (V,A)$ in which the in-neighbourhood as well as the out-neighbourhood of each 
vertex induces a tournament. Such oriented graphs are called 
{\em local tournaments}. Local tournaments are a well-studied class of 
oriented graphs \cite{bangJGT14,huang,mlgx}. Skrien \cite{skrien} proved that 
a connected graph is proper circular-arc graph if and only if it can be oriented
as a local tournament. In fact, every proper circular-arc graph admits a stronger
orientation called a {\em local transitive tournament}. This is an oriented graph
in which the in-neighbourhood as well as the out-neighbourhood of each
vertex induces a transitive tournament. It is proved in \cite{hh} a connected
graph is a proper circular-arc graph if and only if it can be oriented as a local
transitive tournament.      

An oriented graph is {\em acyclic} if it does not contain a directed cycle. Acyclic
local tournaments are a subclass of local transitive tournaments.
The underlying graphs of acyclic local tournaments correspond to a subclass of 
proper circular-arc graphs.
A graph $G = (V,E)$ is a {\em proper interval graph} if there is a family of 
intervals $I_v, v \in V$ such that no interval contains another and for all 
$u, v \in V$, $uv \in E$ if and only if $I_u$ and $I_v$ intersect. A graph is 
a proper interval graph if and only if it has an acyclic local tournament 
orientation, cf. \cite{hh}. 

A {\em partially oriented graph} is a mixed graph $H$ obtained from some graph $G$
by orienting the edges in a subset of the edge set of $G$. The graph $G$ is called
the {\em underlying graph} of $H$. We denote $H$ by $(V,E\cup A)$ where $E$ is
the set of (non-oriented) edges and $A$ is the set of arcs in $H$.
We use $uv$ to denote an edge in $E$ with endvertices $u, v$ and use $(u,v)$ to
denote an arc in $A$ with {\em tail} $u$ and {\em head} $v$.
In either case we say that $u, v$ are {\em adjacent} in $H$.
The partially oriented graph $H$ is {\em connected} if its underlying graph $G$ is.

A class $\cal C$ of graphs is called {\em hereditary} if it is closed under taking 
induced subgraphs, that is, if $G \in \cal C$ and $G'$ is an induced subgraph of 
$G$ then $G' \in \cal C$. Proper interval graphs and proper circular-arc graphs
are examples of hereditary classes of graphs. Hereditary classes of oriented graphs
are defined analogously. We extend this concept to partially oriented graphs.

Let $H = (V,E\cup A)$ and $H' = (V',E'\cup A')$ be partially oriented graphs.
We says that $H$ {\em critically contains} $H'$ (or $H'$ is {\em critically 
contained} in $H$) if $V' \subseteq V$ and for all $u, v \in V'$,
\begin{itemize}
\item $u$ and $v$ are adjacent in $H'$ if and only if they are adjacent in $H$;
\item if $(u,v) \in A'$ then $(u,v) \in A$;
\item if $uv \in E'$, then $uv \in E$, or $(u,v) \in A$, or $(v,u) \in A$.
\end{itemize}.

\vspace{-5mm}

\noindent Equivalently, $H'$ is critically contained $H$ if and only if it is 
obtained from $H$ by possibly deleting vertices, followed by replacing arcs 
with edges.

We note that, in case when $H$ and $H'$ are both graphs or oriented graphs, $H$ 
critically contains $H'$ if and only if $H$ contains $H'$ as an induced subgraph 
or an induced oriented subgraph. We call a class $\cal C$ of partially oriented 
graphs {\em hereditary} if $H \in \cal C$ and $H'$ is critically contained in $H$ 
then $H' \in \cal C$.
 
Let $\cal C$ be a hereditary class of oriented graphs.
The {\em orientation completion problem} for $\cal C$ asks whether a given 
partially oriented graph $H = (V,E\cup A)$ can be completed to an oriented graph 
in $\cal C$ by orienting the edges in $E$. 
The hereditary property of $\cal C$ ensures that if a partially oriented graph $H$ 
can be completed to an oriented graph in $\cal C$ then every partially 
oriented graph that is critically contained in $H$ can also be completed to 
an oriented graph in $\cal C$. Therefore the partially oriented graphs which can be
completed to oriented graphs in $\cal C$ form a hereditary class. 

Orientation completion problems have been studied for several classes of
oriented graphs, including local tournaments, local transitive tournaments, and
acyclic local tournaments, cf. \cite{bhz,huang_locp}. It is proved in \cite{bhz} 
that the orientation completion problem is polynomial time solvable for local 
tournaments and for acyclic local tournaments, but NP-complete for local 
transitive tournaments.

Any hereditary class of graphs admits a characterization by forbidden subgraphs. 
The forbidden subgraphs consists of minimal graphs which do not belong to the class.
This is also the case for a hereditary class of partially oriented graphs and 
in particular for the class of partially oriented graphs which can be completed to 
acyclic local tournaments.   

We call a partially oriented graph $X = (V,E\cup A)$ an {\em obstruction for 
acyclic local tournament orientation completions} (or simply, an {\em obstruction})
if the following three properties hold:

\begin{enumerate}
\item $X$ cannot be completed to an acyclic local tournament;
\item For each $v \in V$, $X-v$ can be completed to an acyclic local tournament;
\item For each $(u,v) \in A$, the partially oriented graph obtained from
           $X$ by replacing $(u,v)$ with the edge $uv$ can be completed to
           an acyclic local tournament.
\end{enumerate}

Thus an obstruction $X$ for acyclic local tournament orientation completions
is a partially oriented graph which cannot be completed to a local tournament and 
is minimal in the sense that if $X'$ is critically contained in $X$ and $X' \neq X$
then $X'$ can be completed to an acyclic local tournament. 
Obstructions for local tournament orientation completions are defined analogously. 
They are the minimal (in the sense of critical containment) partially oriented 
graphs which cannot be completed to local tournaments.

The {\em dual} of an obstruction $X$ is the partially oriented graph obtained from 
$X$ by reversing the arcs in $X$. Clearly, the dual of an obstruction is again 
an obstruction. Obstructions are present in any partially oriented graph that 
cannot be completed to an acyclic local tournament.

\begin{proposition} \label{obstruction}
A partially oriented graph $H$ cannot be completed to an acyclic local tournament 
if and only if it critically contains an obstruction.
\end{proposition}
\pf If $H$ can be completed to an acyclic local tournament, then every partially 
oriented graph critically contained in $H$ can also be completed to an acyclic 
local tournament so $H$ does not critically contain an obstruction. 
On the other hand, 
suppose that $H$ cannot be completed to an acyclic local tournament. 
By deleting vertices and replacing arcs with edges in $H$ as long as the resulting 
partially oriented graph still cannot be completed to an acyclic local tournament 
we obtain an obstruction that is critically contained in $H$.
\qed

We will find in this paper all obstructions for acyclic local tournament 
orientation completions. In particular, we will prove the following:

\begin{theorem} \label{main}
Let $X$ be an obstruction for acyclic local tournament orientations. Then $X$ or 
its dual is a $C_k$ ($k \geq 4$) or one of the graphs in 
Figures~\ref{pigfigure}--\ref{notlt}.
\end{theorem}

In the companion paper \cite{hsu} we will determine all obstructions for local 
tournament orientation completions.

\section{Preliminary results} \label{Preliminary Results}

A {\em straight enumeration} of a graph $G$ is a vertex ordering $\prec$ such that
for all $u \prec v \prec w$, if $uw$ is an edge of $G$, then both $uv$ and $vw$ are
edges. This property is referred to as the {\em umbrella property} of the vertex
ordering.

\begin{theorem} \label{alt} \cite{hh,huang}
The following statements are equivalent for a graph $G$:
\begin{enumerate}
\item $G$ can be completed to an acyclic local tournament;
\item $G$ is a proper interval graph;
\item $G$ has a straight enumeration.
\qed
\end{enumerate}
\end{theorem}

Wegner \cite{wegner} found all minimal graphs which are not proper interval graphs.

\begin{theorem}\label{pig} \cite{wegner}
A graph $G$ is a proper interval graph if and only if it does not contain a $C_k$
($k \geq 4$), a tent, a claw, or a net in Figure~\ref{pigfigure} as an induced 
subgraph. \qed
        \begin{figure}[H]
                \centering
                \captionsetup[subfigure]{labelformat=empty}
                \begin{subfigure}{0.33\textwidth}
                        \centering
                        \begin{tikzpicture}
                        \node[unlabelled] (a) at (1, 1.73)              {};
                        \node[unlabelled] (b) at (0.5, 0.865)   {};
                        \node[unlabelled] (c) at (1.5, 0.865)   {};
                        \node[unlabelled] (d) at (0, 0)                 {};
                        \node[unlabelled] (e) at (1, 0)                 {};
                        \node[unlabelled] (f) at (2, 0)                 {};
                        \draw   (a) edge (b)    (a) edge (c)    (b) edge (c)
                        (b) edge (d)    (b) edge (e)    (c) edge (e)
                        (c) edge (f)    (d) edge (e)    (e) edge (f);
                        \end{tikzpicture}
                        \subcaption{Tent}
                \end{subfigure}%
                \begin{subfigure}{0.33\textwidth}
                        \centering
                        \begin{tikzpicture}
                        \node[unlabelled] (1) at (0:0) {};
                        \node[unlabelled] (2) at (90:1) {};
                        \node[unlabelled] (3) at (210:1) {};
                        \node[unlabelled] (4) at (330:1) {};
                        \draw (1) -- (2) (1) -- (3) (1) -- (4);
                        \end{tikzpicture}
                        \subcaption{Claw}
                \end{subfigure}%
                \begin{subfigure}{0.33\textwidth}
                        \centering
                        \begin{tikzpicture}
                        \node[unlabelled] (1) at (0.5,0.866) {};
                        \node[unlabelled] (3) at (0,0) {};
                        \node[unlabelled] (2) at (1,0) {};
                        \node[unlabelled] (u) at (0.5,1.866) {};
                        \node[unlabelled] (w) at (-1,0) {};
                        \node[unlabelled] (v) at (2,0) {};
                        \draw (1) -- (2) (2) -- (3) (3) -- (1) (1) -- (u) (2) -- (v) (3) -- (w);
                        \end{tikzpicture}
                        \subcaption{Net}
                \end{subfigure}%
                \caption{\label{pigfigure}}
        \end{figure}

\end{theorem}

Theorems \ref{alt} and \ref{pig} imply that $C_{k}$ ($k \geq 4$) and the graphs 
in Figure~\ref{pigfigure} are precisely the obstructions for acyclic local 
tournament orientation completions which contain no arcs. Hence we only need to 
find obstructions that contain arcs. By definition the underlying graph of any 
obstruction that contains arcs is a proper interval graph and hence has a 
straight enumeration.

Let $G = (V,E)$ be a graph and $Z(G) = \{(u,v): uv \in E\}$ be the set of all 
ordered pairs $(u,v)$ such that $uv \in E$. Note that each edge $uv \in E$ 
gives rise to two ordered pairs $(u,v), (v,u)$ in $Z(G)$. Suppose that
$(u, v)$ and $(x, y)$ are two ordered pairs of $Z(G)$. We say $(u,v)$ 
\textit{forces} $(x,y)$ and write $(u,v) \Gamma (x,y)$ if one of the following 
conditions is satisfied:
	
\begin{itemize}
\item $u = x$ and $v = y$;
\item $u = y$, $v \neq x$, and $vx \notin E$;
\item $v = x$, $u \neq y$, and $uy \notin E$.
\end{itemize}
	
\noindent We say that $(u,v)$ \textit{implies} $(x,y)$ and write 
$(u,v) \Gamma^* (x,y)$ if there exists a sequence of pairs 
$(u_1, v_1), (u_2, v_2), \dots, (u_k, v_k) \in Z(G)$ such that 
$$(u, v) = (u_1, v_1) \Gamma (u_2, v_2) \Gamma \dots \Gamma (u_k, v_k) = (x, y).$$
\noindent We will call such a sequence a $\Gamma$-\textit{sequence} from 
$(u, v)$ to $(x, y)$. It is easy to verify that $\Gamma^*$ is an equivalence 
relation on $Z(G)$. 

\begin{proposition}[] \cite{huang} \label{Gamma-class}
Let $G$ be a graph and $D = (V,A)$ be a local tournament orientation of $G$.
Suppose that $(u, v) \Gamma^* (x, y)$ for some $(u, v), (x, y) \in Z(G)$. Then
$(u,v) \in A$ if and only if $(x,y) \in A$.
\qed
\end{proposition}

\begin{proposition}\label{pig_gamma_sequence}
Let $G = (V,E)$ be a proper interval graph and $\prec$ be a straight enumeration of
$G$. Suppose that $(u,v) \Gamma^* (x,y)$. Then $u \prec v$ if and only if 
$x \prec y$.
\end{proposition}
\pf It suffices to show that if $u \prec v$ and $(u,v) \Gamma (x,y)$ then
$x \prec y$. So assume that $(u,v) \Gamma (x,y)$. Then one of the following holds:
\begin{itemize}
\item $u = x$ and $v = y$;
\item $u = y$, $v \neq x$, and $vx \notin E$;
\item $v = x$, $u \neq y$, and $uy \notin E$.
\end{itemize}
Clearly, $x \prec y$ when $u = x$ and $v = y$.
Suppose that $u = y$, $v \neq x$, and $vx \notin E$. If $u \prec x \prec v$,
then it violates the umbrella property because $uv \in E$ but $xv \notin E$.
If $u \prec v \prec x$, then it again violates the umbrella property because
$ux \in E$ but $vx \notin E$. Hence we must have $x \prec u = y$.
The proof for the case when $v = x$, $u \neq y$, and $uy \notin E$ is similar.
\qed

The relation $\Gamma^*$ on $Z(G)$ induces a partition of the edge set of $G$ into
\textit{implication classes} as follows: two edges $uv, xy$ of $G$
are in the same implication class if and only if $(u,v) \Gamma^* (x,y)$ or
$(u,v) \Gamma^* (y,x)$. An implication class is called {\em trivial} if it has
only one edge and {\em non-trivial} otherwise. An edge $uv$ of $G$ is called
{\em balanced} if $N[u] = N[v]$ and {\em unbalanced} otherwise. Clearly,
any balanced edge forms a trivial implication class, and the unique edge in any
trivial implication class is balanced.

A vertex in a graph is {\em universal} if it is adjacent to every other vertex.

\begin{proposition} \cite{huang} \label{pigstructure}
Suppose that $G = (V,E)$ is a connected proper interval graph that is not a
complete graph. Then $\overline{G}$ has a unique non-trivial component $H$.
If $F$ is an implication class of $G$, then $F$ is one of the following types:
\begin{itemize}
\item $F$ is trivial;
\item $F$ consists of all unbalanced edges within $H$;
\item $F$ consists of all edges of $G$ between $H$ and a universal vertex of $G$.
\end{itemize}
In particular, if $G$ contains no universal vertex, then $G$ has a unique
non-trivial implication class.
\qed
\end{proposition}

Let $H$ be a partially oriented graph whose underlying graph $U(H)$ is a proper
interval graph. Suppose that $\prec$ is a straight enumeration of $U(H)$. We call
an arc $(u,v)$ of $H$ \textit{positive} (with respect to $\prec$) if $u \prec v$
and \textit{negative} otherwise. If $H$ does not contain negative arcs, then $H$
can be completed to an acyclic local tournament by replacing all edges of $H$ with
positive arcs. Similarly, if $H$ does not contain positive arcs then it can also
be completed to an acyclic local tournament. 

For convenience we call an arc $(u,v)$ of a partially oriented graph $H$ 
{\em balanced} if $uv$ is balanced in $U(H)$, and {\em unbalanced} otherwise. 
If $H$ contains a directed cycle, then clearly $H$ cannot be completed to 
an acyclic local tournament. On the other hand, when $H$ does not contain a 
directed cycle, whether $H$ can be completed to an acyclic local tournament 
can be recognized from the directions of the unbalanced arcs of $H$ in a straight 
enumeration. The following proposition is a reformulation of a result 
(Corollary 3.3) from \cite{huang}.

\begin{proposition} \cite{huang} \label{unbalanced}
Let $H$ be a partially oriented graph such that $U(H)$ is a proper interval graph
and let $\prec$ be a straight enumeration of $U(H)$. Suppose $H$ does not contain
a directed cycle. Then $H$ cannot be completed to an acyclic local tournament if 
and only if it contains two unbalanced arcs, one positive and one negative with 
respect to $\prec$.
\qed
\end{proposition}

Let $v$ be a vertex and $(x,y)$ be an arc in a partially oriented graph $H$ where 
$v \notin \{x,y\}$. We call $v$ the {\em $(x,y)$-balancing} vertex if $v$ is 
the only vertex adjacent to exactly one of $x,y$, that is, $(x,y)$ is not balanced
in $H$ and is balanced in $H-v$. When $v$ is an $(x,y)$-balancing vertex, we call 
the set $\{v,x,y\}$ an {\em arc-balancing triple} and in the case that the arc 
$(x,y)$ does not need to be specified, we simply call $v$ an {\em arc-balancing} 
vertex.

\section{Obstructions}

It remains to determine the obstructions for acyclic local tournament orientation
completions that contain arcs. By Theorem \ref{alt} their underlying graphs 
are proper interval graphs and hence have straight enumerations.
Of these obstructions some cannot even be completed to local tournaments and the 
rest can be completed to local tournaments but not to acyclic local tournaments. 
We will treat these two types of obstructions separately
(See Theorems \ref{cancompletetheorem} and \ref{cannotcompletetheorem} below). 
Note that any obstruction for acyclic local tournament orientation completions 
that cannot be complete to a local tournament is an obstruction for local 
tournament orientation completions by definition.    

\begin{lemma}\label{cancompletelemma}
Let $X$ be an obstruction for acyclic local tournament orientation completions that
contains arcs but no directed cycle. Then any vertex not incident with an arc is 
an arc-balancing vertex. Moreover, if $X$ can be completed to a local tournament
then there exists a universal vertex incident with exactly one arc of $X$.
\end{lemma}
\pf Let $\prec$ be a straight enumeration of $U(X)$. By Proposition \ref{unbalanced}
$X$ contains two unbalanced arcs $(a,b), (c,d)$, one positive and one negative 
with respect to $\prec$. Since $X$ is an obstruction for acyclic local 
tournament orientation completions, $X-v$ can be completed to an acyclic local
tournament for any vertex $v$. If $v \notin \{a,b,c,d\}$, then it follows from
Proposition \ref{unbalanced} that at least one of $(a,b), (c,d)$ is balanced in 
$X-v$, which means that $v$ is an arc-balancing vertex. 

Suppose that $X$ can be completed to a local tournament. By Propositions
\ref{Gamma-class} and \ref{pig_gamma_sequence}, $ab$ and $cd$ belong to different 
implication classes of $U(X)$. Hence, by Proposition \ref{pigstructure}, one of 
$ab, cd$ is an edge of $U(X)$ between the unique non-trivial component of 
$\overline{U(X)}$ and a universal vertex of $U(X)$. That is, there exists a 
universal vertex incident with exactly one of $(a,b), (c,d)$.
\qed

\begin{lemma} \label{dicycle}
Let $X$ be an obstruction for acyclic local tournament orientation completions.
If $X$ contains a directed cycle, then $X$ is Figure~\ref{cancomplete}(viii).
\end{lemma}
\pf Let $C: u_1u_2 \dots u_n$ be a shortest directed cycle in $X$. The cycle $C$ 
must contain all vertices of $X$ as otherwise there is a vertex $v$ not on $C$ 
such that $X - v$ cannot be completed to an acyclic local tournament, 
a contradiction to the assumption that $X$ is an obstruction. 
Since $C$ is a shortest directed cycle containing all vertices of $X$, $X$ does 
not contain arcs other than those of $C$. Since $X$ is an obstruction,
the partially oriented graph obtained from $X$ by replacing $(u_n,u_1)$ with 
the edge $u_nu_1$ can be completed to an acyclic local tournament $D$. 
Now $D$ is acyclic and contains the directed path $u_1u_2 \dots u_n$. Hence
for all adjacent vertices $u_i, u_j$ with $i < j$, $(u_i,u_j)$ is an arc in $D$. 
In particular, $(u_1,u_n)$ is an arc in $D$.
Suppose $D$ contains a pair of nonadjacent vertices. Let $u_i, u_j$ be such a pair
with $j-i$ being maximum. Then $i>1$ or $j<n$. In the case when $i> 1$,
the choice of $u_i,u_j$ implies that $(u_{i-1},u_j)$ is an arc. Now
$(u_{i-1},u_i)$ and $(u_{i-1},u_j)$ are arcs in $D$. Since $D$ is a local
tournament $u_i, u_j$ must be adjacent, which is a contradiction.
A similar proof applies to the case $j<n$ and leads to a contradiction.
Therefore the vertices in $D$ are pairwise adjacent and we see that
$X$ is Figure~\ref{cancomplete}(viii). 
\qed 

\begin{theorem}\label{cancompletetheorem}
Let $X$ be an obstruction for acyclic local tournament orientation completions 
that contains arcs. Suppose that $X$ can be completed to a local tournament. 
Then $X$ or its dual is one of the graphs in Figure~\ref{cancomplete}.
	
	\begin{figure}[H]
		\centering
		\captionsetup[subfigure]{labelformat=empty}
		\begin{subfigure}[b]{0.33\textwidth}
			\centering
			\begin{tikzpicture}
			\node[unlabelled] (1) at (0,0) {};
			\node[unlabelled] (2) at (1,0) {};
			\node[unlabelled] (3) at (2,0) {};
			\node[unlabelled] (4) at (3,0) {};
			\draw (2) edge (3) (1) edge[bend left = 45] (3) (2) edge[bend left = 45] (4);
			\draw[->-=0.5] (1) to (2);
			\draw[->-=0.5] (4) to (3);
			
			\end{tikzpicture}
			\subcaption{(i)}
		\end{subfigure}%
		\begin{subfigure}[b]{0.33\textwidth}
			\centering
			\begin{tikzpicture}
			\node[unlabelled] (1) at (0,0) {};
			\node[unlabelled] (2) at (1,0) {};
			\node[unlabelled] (3) at (2,0) {};
			\node[unlabelled] (4) at (3,0) {};
			\draw (2) -- (3) (3) -- (4);
			\draw (2) edge[bend left = 45] (4);
			\draw[->-=0.5] (1) to (2); 
			\draw[->-=0.5] (3) to[bend right=45] (1);
			\end{tikzpicture}
			\subcaption{(ii)}
		\end{subfigure}%
		\begin{subfigure}[b]{0.33\textwidth}
			\centering
			\begin{tikzpicture}
			\node[unlabelled] (1) at (0,0) {};
			\node[unlabelled] (2) at (1,0) {};
			\node[unlabelled] (3) at (2,0) {};
			\node[unlabelled] (4) at (3,0) {};
			\node[unlabelled] (5) at (4,0) {};
			\draw (3) -- (4) (4) -- (5);
			\draw (1) edge[bend left = 45] (3) (2) edge[bend left = 45] (4) (3) edge[bend left = 45] (5);
			\draw[->-=0.5] (1) to (2);
			\draw[->-=0.5] (3) to (2);
			\end{tikzpicture}
			\subcaption{(iii)}
		\end{subfigure}%
		
		\begin{subfigure}[b]{0.33\textwidth}
			\centering
			\begin{tikzpicture}
			\node[unlabelled] (1) at (0,0) {};
			\node[unlabelled] (2) at (1,0) {};
			\node[unlabelled] (3) at (2,0) {};
			\node[unlabelled] (4) at (3,0) {};
			\node[unlabelled] (5) at (4,0) {};
			\draw (5) -- (4) (2) -- (3);
			\draw (3) edge (4) (2) edge[bend left = 45] (4) (3) edge[bend left = 45] (5) (1) edge[bend left=45] (3);
			\draw[->-=0.5] (5) to[bend right = 45] (3);
			\draw[->-=0.5] (1) to (2);
			\end{tikzpicture}
			\subcaption{(iv)}
		\end{subfigure}%
		\begin{subfigure}[b]{0.33\textwidth}
			\centering
			\begin{tikzpicture}
			\node[unlabelled] (1) at (0,0) {};
			\node[unlabelled] (2) at (1,0) {};
			\node[unlabelled] (3) at (2,0) {};
			\node[unlabelled] (4) at (3,0) {};
			\node[unlabelled] (5) at (4,0) {};
			\draw[->-=0.5] (1) to[bend left = 45] (3);
			\draw[->-=0.5] (4) to[bend right = 45] (2);
			\draw (1) edge (2) (2) edge (3) (3) edge (4) (4) edge (5) (3) edge[bend left = 45] (5) (1) edge[bend left = 55] (4);
			\end{tikzpicture}
			\subcaption{(v)}
		\end{subfigure}%
		\begin{subfigure}[b]{0.33\textwidth}
			\centering
			\begin{tikzpicture}
			\node[unlabelled] (1) at (0,0) {};
			\node[unlabelled] (2) at (1,0) {};
			\node[unlabelled] (3) at (2,0) {};
			\node[unlabelled] (4) at (3,0) {};
			\node[unlabelled] (5) at (4,0) {};
			\node[unlabelled] (6) at (5,0) {};
			\draw (6) -- (5) (2) -- (3) (4) -- (5);
			\draw (1) edge[bend left = 45] (3) (2) edge[bend left = 45] (4) (3) edge[bend left = 45] (5) (4) edge[bend left = 45] (6) (2) edge[bend left = 55] (5) (1) edge[bend left = 55] (4);
			\draw[->-=0.5] (4) to (3);
			\draw[->-=0.5] (1) to (2);
			\end{tikzpicture}
			\subcaption{(vi)}
		\end{subfigure}%
		
		\begin{subfigure}[b]{0.4\textwidth}
			\centering
			\begin{tikzpicture}
			\node[unlabelled] (1) at (0,0) {};
			\node[unlabelled] (2) at (1,0) {};
			\node[unlabelled] (3) at (2,0) {};
			\node[unlabelled] (4) at (3,0) {};
			\node[unlabelled] (5) at (4,0) {};
			\node[unlabelled] (6) at (5,0) {};
			\draw (1) edge (2) (3) edge (4) (5) edge (6) (1) edge[bend left = 45] (3) (2) edge[bend left = 45] (4) (3) edge[bend left = 45] (5) (4) edge[bend left = 45] (6) (1) edge [bend left = 55] (4) (2) edge [bend left = 55] (5) (3) edge [bend left = 55] (6);
			\draw[->-=0.5] (2) to (3);
			\draw[->-=0.5] (5) to (4);
			\end{tikzpicture}
			\subcaption{(vii)}
		\end{subfigure}%
		\begin{subfigure}[b]{0.4\textwidth}
			\centering
			\begin{tikzpicture}
			\node[unlabelled] (1) at (0,0) {};
			\node[unlabelled] (2) at (1,0) {};
			\node (3) at (2,0) {$\dots$};
			\node[unlabelled] (4) at (3,0) {};
			\node[unlabelled] (5) at (4,0) {};
			\draw (1) edge[bend left = 45] (3) (2) edge[bend left = 45] (4) (3) edge[bend left = 45] (5) (1) edge[bend left = 55] (4) (2) edge[bend left = 55] (5);
			\draw[->-=0.5] (1) to (2);
			\draw[->-=0.5] (2) to (3);
			\draw[->-=0.5] (3) to (4);
			\draw[->-=0.5] (4) to (5);
			\draw[->-=0.5] (5) to[bend right = 65] (1);
			\end{tikzpicture}
			\subcaption{(viii)}
		\end{subfigure}%
		
\caption{Obstructions that contain arcs and can be completed to local tournaments. \label{cancomplete}}
	\end{figure}
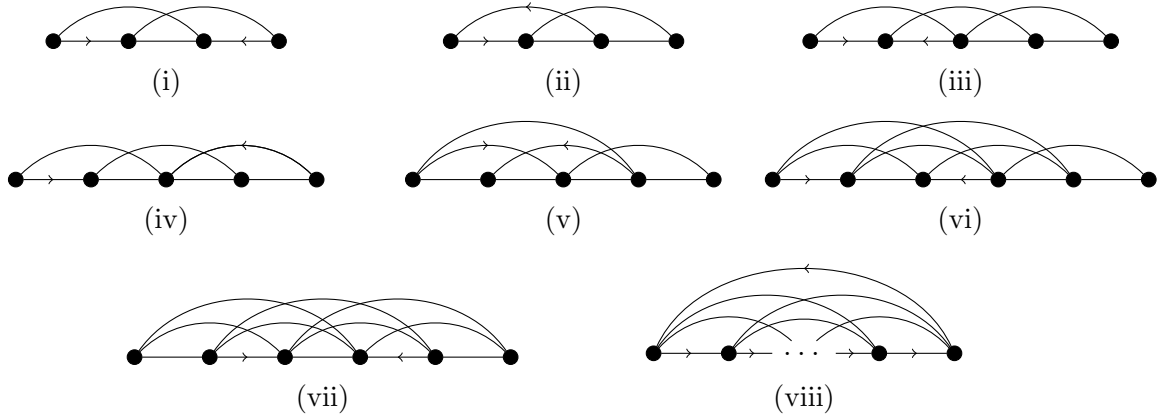
\end{theorem}
\pf It is easy to verify that each graph in Figures~\ref{cancomplete} is 
an obstruction for acyclic local tournament orientation completions and can be 
completed to a local tournament. Hence it suffices to show that $X$ is one of them.
If $X$ contains a directed cycle, then it is Figure~\ref{cancomplete}(viii) by 
Lemma \ref{dicycle}. So we may assume that $X$ does not contain a directed cycle. 

Let $\prec:v_1, v_2, \dots, v_n$ be a straight enumeration of $U(X)$. 
By Proposition \ref{unbalanced}, $X$ contains exactly two unbalanced arcs, one 
positive and one negative with respect to $\prec$. By Lemma \ref{cancompletelemma}, 
there exists a universal vertex $v_c$ incident with exactly one arc of $X$. 
Since $X$ contains unbalanced arcs, $U(X)$ is not a complete graph. This implies
$v_1v_n$ is not an edge of $U(X)$ and in particular, $c \notin \{1, n\}$.

Clearly, $n \geq 3$. If $n=3$, then $X$ or its dual must contain the arcs 
$(v_1,v_2)$ and $(v_3, v_2)$, contradicting the fact that $X$ can be completed to 
a local tournament. 
Suppose $n=4$. Without loss of generality, assume $c = 2$. So, 
$v_2v_4$ is an edge of $U(X)$. If $v_1v_3$ is not an edge of $U(X)$, then the only 
unbalanced edges of $U(X)$ are those incident with $v_2$, so both unbalanced arcs 
of $X$ are incident with $v_2$, contradicting the choice of $v_c$. Hence, $v_1v_3$ 
is an edge. Since both arcs of $X$ are unbalanced, there is no arc between $v_2$ 
and $v_3$. It is now easy to see that $X$ or its dual is one of 
Figure~\ref{cancomplete}(i) or (ii).

Suppose $n=5$. We claim that $v_1v_3, v_2v_4, v_3v_5$ are edges of $U(X)$. 
Indeed, if any of $v_1v_3, v_2v_4, v_3v_5$ is not an edge of of $U(X)$, then 
there is an universal vertex incident with both (unbalanced) arcs of $X$, 
a contradiction to Lemma \ref{cancompletelemma}.  
Suppose neither of $v_1v_4, v_2v_5$ is an edge of $U(X)$. Then, $v_3$ is 
the unique universal vertex, so $c = 3$. 
If the arc not incident with $v_3$ is between $v_2$ and $v_4$, then 
$X$ or its dual critically contains Figure 2(i) or (ii), a contradiction to
the minimality of $X$. Hence the arc not incident with $v_3$ is either between 
$v_1$ and $v_2$ or between $v_4$ and $v_5$. We may assume without loss of 
generality $(v_1, v_2)$ is an arc. If $v_5$ is not incident with 
an arc, then it is an arc-balancing vertex by Lemma \ref{cancompletelemma}. Clearly,
$v_5$ must balance the arc $(v_3, v_2)$, so $X$ is Figure~\ref{cancomplete}(iii). 
Otherwise $v_5$ is incident with an arc and $X$ is Figure~\ref{cancomplete}(iv). 
Suppose instead that $v_1v_4$ or $v_2v_5$ is an edge of $U(X)$. Without loss of 
generality, assume $v_1v_4$ is an edge. 
If $v_2v_5$ is also an edge, then each of $v_2, v_3, v_4$ is a universal vertex 
and hence is not an arc-balancing vertex. By Lemma \ref{cancompletelemma}, 
each of $v_2, v_3, v_4$ is incident with an arc, so there is an arc with both 
endvertices among $v_2, v_3, v_4$, contradicting the fact that both arcs are 
unbalanced. So, $v_2v_5$ is not an edge. Any arc incident with $v_5$ does not have 
an arc-balancing vertex because there are two vertices adjacent to exactly one 
endvertex of such an arc. If $v_1$ or $v_2$ is an arc-balancing vertex, then 
it balances an arc incident with $v_5$, so neither $v_1$ nor $v_2$ is 
an arc-balancing vertex. By Lemma \ref{cancompletelemma}, both $v_1$ and $v_2$ are 
incident with arcs. Similarly, neither $v_3$ nor $v_4$ are arc-balancing vertices 
because they are universal, so they are both incident with arcs. Since both arcs 
of $X$ are unbalanced, $X$ or its dual must be Figure~\ref{cancomplete}(v).

Suppose $n \geq 6$. Since $X$ contains exactly two arcs, it contains 
at most two arc-balancing vertices. Since any vertex not incident with an arc is 
an arc-balancing vertex by Lemma \ref{cancompletelemma}, $X$ contains at most two 
vertices not incident with arcs. In particular, $n = 6$ and $X$ contains two 
disjoint arc-balancing triples. We show that neither $v_2$ nor $v_5$ is universal. 
Assume $v_2$ is universal. Since $X$ contains two disjoint arc-balancing triples, 
one of them contains only vertices succeeding $v_1$. Since $v_2v_6$ is an edge of 
$U(X)$, this arc-balancing triple induces a clique in $U(X)$ by the umbrella 
property, a contradiction. Hence, neither $v_2$ nor $v_5$ is universal by symmetry.
So, $c \in \{3,4\}$. Assume $c = 4$ without loss of generality. Let $v_k$ be 
the arc-balancing vertex for the arc incident with $v_4$ and $v_j$ be the other 
endvertex. Then, $v_k$ is the unique vertex adjacent to $v_4$ and not $v_j$, so 
$v_j$ is adjacent to every vertex except for $v_k$. It follows from the straight 
enumeration that $k \in \{1,6\}$.

Suppose $k = 6$. If $v_j = v_1$, then $v_1$ is adjacent to $v_5$, contradicting
the fact that $v_5$ is not a universal vertex, so $v_j \neq v_1$. Since $v_j$ is 
not adjacent to $v_6$, we have either $v_j = v_3$ or $v_j = v_2$. First suppose 
$v_j = v_3$. Without loss of generality, assume $v_6$ is a $(v_4, v_3)$-balancing 
vertex. Since $X$ contains two disjoint arc-balancing triples, $\{v_1, v_2, v_5\}$ 
is an arc-balancing triple. If $v_1$ balances an arc between $v_2$ and $v_5$, then 
$v_6$ is adjacent to both $v_2$ and $v_5$, contradicting the fact that $v_2$ is 
not a universal vertex. Clearly, $v_2$ cannot balance an arc between $v_1$ and 
$v_5$ by the straight enumeration. So, $v_5$ is a $(v_1, v_2)$-balancing vertex 
and thus $X$ is Figure~\ref{cancomplete}(vi). On the other hand, suppose 
$v_j = v_2$. Without loss of generality, assume $v_6$ is a $(v_4, v_2)$-balancing 
vertex. By a similar argument as above, $\{v_1, v_3, v_5\}$ is an arc-balancing 
triple. Clearly, $v_3$ cannot be the arc-balancing vertex by the straight 
enumeration. If $v_1$ is a $(v_3, v_5)$-balancing vertex, then the dual of $X$ is 
Figure~\ref{cancomplete}(vii). Otherwise $v_5$ is a $(v_1, v_3)$-balancing vertex 
and $X$ is Figure~\ref{cancomplete}(vi).

Suppose now $k = 1$. If $v_j = v_6$, then $v_j$ is adjacent to $v_2$, so $v_2$ is 
a universal vertex, a contradiction. Hence, $v_j \neq v_6$. Since $v_j$ is not 
adjacent to $v_k$, we have $v_j = v_5$. Without loss of generality, assume 
$v_1$ is a $(v_5, v_4)$-balancing vertex. By a similar argument as above, 
$\{v_2, v_3, v_6\}$ is an arc-balancing triple. If $v_2$ balances an arc between 
$v_3$ and $v_6$, then $v_6$ must be adjacent to $v_1$, a contradiction. Clearly, 
$v_3$ cannot balance an arc between $v_2$ and $v_6$ by the straight enumeration. 
Hence, $v_6$ is a $(v_2, v_3)$-balancing vertex. It is now easy to see that $X$ is 
Figure~\ref{cancomplete}(vii).
\qed

\begin{figure}[H]
	\centering
	\captionsetup[subfigure]{labelformat=empty}
	\begin{subfigure}[b]{.33\textwidth}
		\centering
		\begin{tikzpicture}
			\node[unlabelled] (a1) at (3,0) {};
			\node[unlabelled] (a2) at (2,0) {};
			\node[unlabelled] (x1) at (1,0) {};
			\node[unlabelled] (b1) at (0,0) {};
			\draw	(a1) -- (a2)	(x1) -- (b1);
			\draw[->-=.5]	(a1) to [bend right = 45] (x1);
			\draw[->-=.5]	(x1) to (a2);
		\end{tikzpicture}
		\subcaption{(i)}
	\end{subfigure}%
	\begin{subfigure}[b]{.33\textwidth}
		\centering
		\begin{tikzpicture}
			\node[unlabelled]	(1) at (0,0)	{};
			\node[unlabelled]	(2) at (1,0)	{};
			\node[unlabelled]	(3) at (2,0)	{};
			\node[unlabelled]	(4) at (3,0)	{};
			\node[unlabelled]	(5) at (4,0)	{};
			\draw	(1) -- (2)	(4) -- (5)
			(2) edge[bend left = 45] (4);
			\draw[->-=.5]	(2) to (3);
			\draw[->-=.5]	(4) to (3);
		\end{tikzpicture}
		\subcaption{(ii)}
	\end{subfigure}%
	\begin{subfigure}[b]{.33\textwidth}
		\centering
		\begin{tikzpicture}
			\node[unlabelled]	(1) at (0,0)	{};
			\node[unlabelled]	(2) at (1,0)	{};
			\node[unlabelled]	(3) at (2,0)	{};
			\node[unlabelled]	(4) at (3,0)	{};
			\node[unlabelled]	(5) at (4,0)	{};
			\draw	(1) edge (2) (4) edge (5)
			(1) edge[bend left = 45] (3)
			(2) edge[bend left = 45] (4)
			(3) edge[bend left = 45] (5);
			\draw[->-=.5]	(2) to (3);
			\draw[->-=.5]	(4) to (3);
		\end{tikzpicture}
		\subcaption{(iii)}
	\end{subfigure}%

	\begin{subfigure}[b]{0.5\textwidth}
		\centering
		\begin{tikzpicture}
			\node[unlabelled]	(1) at (0,0)	{};
			\node[unlabelled]	(2) at (1,0)	{};
			\node[unlabelled]	(3) at (2,0)	{};
			\node[unlabelled]	(4) at (3,0)	{};
			\node[unlabelled]	(5) at (4,0)	{};
			\node[unlabelled]	(6) at (5,0)	{};
			\draw	(1) -- (2)	(5) -- (6)	(3) -- (4)
			(2) edge[bend left = 45] (4)
			(3) edge[bend left = 45] (5)
			(2) edge[bend left = 55] (5);
			\draw[->-=.5]	(2) to (3);
			\draw[->-=.5]	(5) to (4);
		\end{tikzpicture}
		\subcaption{(iv)}
	\end{subfigure}%
	\begin{subfigure}[b]{0.5\textwidth}
		\centering
		\begin{tikzpicture}
			\node[unlabelled]	(0) at (-1,0)	{};
			\node[unlabelled]	(1) at (0,0)	{};
			\node[unlabelled]	(2) at (1,0)	{};
			\node[unlabelled]	(3) at (2,0)	{};
			\node				(4) at (3,0)	{$\dots$};
			\node[unlabelled]	(6) at (4,0)	{};
			\node[unlabelled]	(7) at (5,0)	{};
			\draw	(0)--(1)	(2)--(3)	(3)--(4)	(4)--(6);
			\draw	(1) edge[bend left=45]	(3);
			\draw[->-=.5]	(1) to (2);
			\draw[->-=.5]	(7) to (6);
		\end{tikzpicture}
		\subcaption{(v)}
	\end{subfigure}%

	\begin{subfigure}[b]{.5\textwidth}
		\centering
		\begin{tikzpicture}
			\node[unlabelled]	(1) at (0,0)	{};
			\node[unlabelled]	(2) at (1,0)	{};
			\node				(4) at (2,0)	{$\dots$};
			\node[unlabelled]	(6) at (3,0)	{};
			\node[unlabelled]	(7) at (4,0)	{};
			\draw	(2)--(4)	(4)--(6); 
			\draw[->-=.5]	(1) to (2);
			\draw[->-=.5]	(7) to (6);
		\end{tikzpicture}
		\subcaption{(vi)}
	\end{subfigure}%
	\begin{subfigure}[b]{0.5\textwidth}
		\centering
		\begin{tikzpicture}
			\node[unlabelled]	(1) at (0,0)	{};
			\node[unlabelled]	(2) at (1,0)	{};
			\node[unlabelled]	(3) at (2,0)	{};
			\node[unlabelled]	(4) at (3,0)	{};
			\node[unlabelled]	(5) at (4,0)	{};
			\node[unlabelled]	(6) at (5,0)	{};
			\draw (2) edge (3) (3) edge (4) (4) edge (5) (1) edge[bend left = 45] (3) (2) edge[bend left = 45] (4) (3) edge[bend left = 45] (5) (4) edge[bend left = 45] (6);
			\draw[->-=.5]	(1) to (2);
			\draw[->-=.5]	(6) to (5);
		\end{tikzpicture}
		\subcaption{(vii)}
	\end{subfigure}

	\begin{subfigure}[b]{\textwidth}
		\centering
		\begin{tikzpicture}
			\node[unlabelled]	(0) at (-1,0)	{};
			\node[unlabelled]	(1) at (0,0)	{};
			\node[unlabelled]	(2) at (1,0)	{};
			\node[unlabelled]	(3) at (2,0)	{};
			\node				(4) at (3,0)	{$\dots$};
			\node[unlabelled]	(5) at (4,0)	{};
			\node[unlabelled]	(6) at (5,0)	{};
			\node[unlabelled]	(7) at (6,0)	{};
			\node[unlabelled]	(8) at (7,0)	{};
			\draw	(0)--(1)	(2)--(3)	(3)--(4)	(4)--(5)
			(5)--(6)	(7)--(8);
			\draw	(1) edge[bend left=45]	(3)
			(5) edge[bend left=45]	(7);
			\draw[->-=.5]	(1) to (2);
			\draw[->-=.5]	(7) to (6);
		\end{tikzpicture}
		\subcaption{(viii)}
	\end{subfigure}%
\caption{Obstructions that contain arcs and cannot be completed to 
local tournaments. \label{notlt}}
\end{figure}
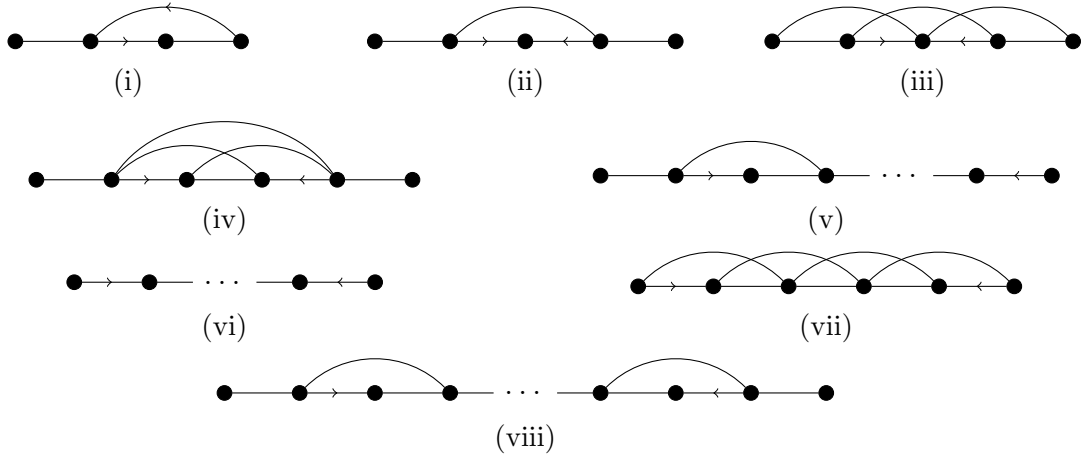

\begin{theorem}\label{cannotcompletetheorem}
Let $X$ be an obstruction for acyclic local tournament orientation completions that
contains arcs. Suppose that $X$ cannot be completed to a local tournament. Then $X$
or its dual is one of the graphs in Figure~\ref{notlt}.
\end{theorem}

\pf It is easy to verify that each graph in Figures~\ref{notlt} is an obstruction 
for acyclic local tournament orientation completions and cannot be completed to 
a local tournament. Hence it suffices to show that $X$ contains one of them as 
an induced subgraph. Since $X$ is an obstruction for acyclic local tournament
orientation completions and cannot be completed to a local tournament, it does not
contain any graph in Figure~\ref{cancomplete} as an induced subgraph. 
In particular, $X$ is not Figure~\ref{cancomplete}(viii) and hence by Lemma 
\ref{dicycle} $X$ does not contain a directed cycle. Therefore by 
Proposition \ref{unbalanced}, $X$ contains exactly two unbalanced arcs, one
positive and one negative with respect to any straight enumeration of $U(X)$.

Let $\prec: v_1, v_2, \dots, v_n$ be a straight enumeration of $U(X)$ and let 
$(v_a,v_b), (v_c,v_d)$ be the two arcs in $X$ where $a < b$ and $c > d$. 
Consider first the case when the two arcs share an endvertex.
Suppose that $v_a = v_d$ is the shared endvertex. By considering the dual of $X$ if 
necessary we assume $b < c$. Then we have $a = d < b < c$ and the umbrella 
property of $\prec$ implies $v_bv_c$ is an edge of $X$. Since $X$ does not contain 
Figure~\ref{cancomplete}(ii), any vertex $v_j$ with $j>c$ adjacent to $v_b$ is 
adjacent to $v_a$. This together with the umbrella property of $\prec$ imply that 
any vertex adjacent to adjacent to $v_b$ is adjacent to $v_a$. The arc $(v_a,v_b)$ 
is unbalanced so there is a vertex $v_i$ adjacent to $v_a$ but not to $v_b$. 
Clearly, we must have $i < a$ and thus the subgraph of $X$ induced by 
$v_i, v_a, v_b, c_c$ is Figures~\ref{notlt}(i). The case when $v_b=v_c$ is
the shared endvertex can be treated analogously. Suppose that $v_b=v_d$ is the 
shared endvertex. If $v_av_c$ is not an edge of $X$, then $v_a, v_b, v_c$ induce 
a Figure 3(vi). So assume $v_av_c$ is an edge.
Since both arcs are unbalanced, for each of them there exists a vertex 
adjacent to exactly one of the two endvertices. Suppose there is a vertex $v_i$
adjacent to $v_a$ but not to $v_b$. Clearly, $i < a$. If there is a vertex $v_j$
adjacent to $v_c$ but not to $v_d = v_b$, then $j>c$ and the subgraph of $X$ 
induced by $v_i, v_a, v_b, v_c, v_j$ is Figure~\ref{notlt}(ii). If there is 
a vertex $v_k$ adjacent to $v_d = v_b$ but not to $v_c$, then $i<k<a$ and the 
subgraph of $X$ induced by $v_i,v_k,v_a,v_b,v_c$ is Figure~\ref{cancomplete}(iii),
a contradiction. Thus we may assume that any vertex adjacent to $v_a$ except $v_b$
is adjacent to $v_b$. Hence there is a vertex adjacent to $v_b$ but not to $v_a$ and
let $v_r$ be such a vertex. If there is a vertex $v_{\ell}$ adjacent to $v_d=v_b$ 
but not to $v_c$, then ${\ell}<a<b<c<r$ and the subgraph of $X$ induced by 
$v_\ell,v_a,v_b,v_c,v_r$ is Figure~\ref{notlt}(iii). If there is a vertex $v_q$
adjacent to $v_c$ but not to $v_d=v_b$, then $a<b<c<r<q$ and the subgraph of $X$ 
induced by $v_{\ell},v_a,v_b,v_c,v_r$ is Figure~\ref{cancomplete}(iii),
a contradiction. The proof for the case when $v_a=v_c$ is the same by considering 
the dual of $X$. Therefore we may further assume the endvertices of the two arcs 
are pairwise distinct.

Suppose that the endvertices of the two arcs are pairwise adjacent. 
Let $v_i$ be a vertex adjacent to exactly one of $v_a,v_b$ and $v_j$ be a vertex
adjacent to exactly one of $v_c,v_d$. Suppose first that $v_i$ is adjacent to 
$v_a$ but not to $v_b$ and $v_j$ is adjacent to $v_d$ but not to $v_c$. Clearly, 
$\mbox{max}\{i,j\}<\mbox{min}\{a,d\}$. The umbrella property of $\prec$ implies 
$v_iv_j$ is an edge of $X$. Thus $v_iv_av_dv_j$ is a $C_4$ in $U(X)$ which cannot 
be induced. So $v_iv_d$ or $v_jv_a$ is an edge of $X$. By symmetry assume 
$v_iv_d$ is an edge. If $v_iv_c$ is not an edge of $X$ then the subgraph induced 
by $v_i,v_a,v_b,v_c,v_d$ is Figure~\ref{cancomplete}(v), a contradiction. So
$v_iv_c$ is an edge, which implies $c < b$. Since $v_j$ is not adjacent to $v_c$
and $c<b$, $v_j$ is not adjacent to $v_b$. If $v_j$ is not adjacent to $v_a$, then 
the subgraph of $X$ induced by $v_i,v_j,v_a,v_b,v_c,v_d$ is
Figure~\ref{cancomplete}(vi), a contradiction. If $v_j$ is adjacent to $v_a$, 
then the subgraph of $X$ induced by $v_j,v_a,v_b,v_c,v_d$ is 
Figure~\ref{cancomplete}(v), a contradiction. 
Suppose now that $v_i$ is adjacent to $v_a$ but not to $v_b$ and $v_j$ is adjacent 
to $v_c$ but not to $v_d$. (Note that the other two cases are symmetric.) 
If $v_i$ is adjacent to neither of $v_c,v_d$ and $v_j$ is adjacent to neither of
$v_a,v_b$, then the subgraph induced by $v_i,v_j,v_a,v_b,v_c,v_d$ is
Figure~\ref{notlt}(iv). If $v_i$ is adjacent to exactly one of $v_c,v_d$, then it 
is adjacent to $v_d$, in which case the subgraph induced by $v_i,v_a,v_b,v_c,v_d$ 
is Figure~\ref{cancomplete}(v), a contradiction. So $v_i$ is adjacent to both
$v_c,v_d$. This implies $c<b$ because $v_iv_b$ is not an edge of $X$. Thus
$v_jv_b$ is an edge following the umbrella property. If $v_j$ is not adjacent to
$v_a$ then the subgraph of $X$ induced by $v_j,v_a,v_b,v_c,v_d$ is 
Figure~\ref{cancomplete}(v), a contradiction. If $v_j$ is adjacent to $v_a$,
then the subgraph of $X$ induced by $v_i,v_j,v_a,v_b,v_c,v_d$ is
Figure~\ref{cancomplete}(vii), a contradiction.

Suppose that the endvertices of the two arcs are not all pairwise adjacent. Without
loss of generality assume $a<d$. Then we must have $b<c$ and in particular $v_av_c$ is not an edge of $X$. Since $X$ does not contain Figure~\ref{cancomplete}(i) as
an induced subgraph, we must have $b<d$ and at least one of $v_av_d$ and $v_bv_c$
is not an edge of $X$. By symmetry we assume $v_av_d$ is not an edge of $X$. 
If $a<b-1$ then $v_{a+1}$ is an arc-balancing vertex according to 
Lemma~\ref{cancompletelemma}. Clearly, $v_{a+1}$ does not balance $(v_a,v_b)$ 
because it is adjacent to both $v_a,v_b$ so it balances $(v_c,v_d)$.
It follows that $v_{a+1}v_d$ is an edge of $X$, which implies
$v_bv_d$ is also an edge of $X$. Since $v_{a+1}$ is the unique vertex adjacent to
exactly one of the endvertices of $(v_c,v_d)$, $v_bv_c$ must be an edge of $X$. 
We see now that the subgraph of $X$ induced by $v_a,v_{a+1},v_b,v_c,v_d$ is
Figure~\ref{cancomplete}(iv), a contradiction. Hence $v_a,v_b$ are consecutive in
$\prec$. Similarly, $v_c,v_d$ are consecutive in $\prec$. If $v_bv_c$ is an edge of
$X$, then any vertex adjacent to $v_d$ except $v_c$ is adjacent to $v_c$. So
there must be a vertex $v_j$ adjacent to $v_c$ but not to $v_d$. The subgraph of 
$X$ induced by $v_a,v_b,v_c,v_d,v_j$ is a graph in Figure~\ref{notlt}(v). 
So we may assume $v_bv_c$ is not an edge of $X$. If $v_bv_d$ is an edge of $X$, 
then the subgraph of $X$ induced by $v_a,v_b,v_c,v_d$ is a graph in 
Figure~\ref{notlt}(vi). So we may further assume $v_bv_d$ is not an edge of $X$. 

Let $v_k$ be the neighbour of $v_b$ having the largest subscript $k$ and let
$v_{\ell}$ be the neighbour of $v_d$ having the least subscript. Clearly,
$b < k < d$ and $b < \ell < d$. Suppose neither $v_av_k$ nor $v_{\ell}v_c$ is an 
edge of $X$. Consider first the case when $\ell < k$. If $v_av_{\ell}$ and 
$v_kv_c$ are both edges of $X$, then the subgraph of $X$ induced by 
$v_a,v_b,v_{\ell},v_k,v_c,v_d$ is Figure~\ref{notlt}(vii). If $v_av_{\ell}$ is not 
an edge of $X$, then the subgraph of $X$ induced by $v_a,v_b,v_{\ell},v_c,v_d$ 
is a graph in Figure~\ref{notlt}(vi). Similarly, if $v_kv_c$ is not an edge of $X$, 
then the subgraph of $X$ induced by $v_a,v_b,v_k,v_c,v_d$ is a graph in 
Figure~\ref{notlt}(vi). When $k \leq \ell$, the subgraph of $X$ induced by 
$v_a,v_b,v_c,v_d$ together with the vertices in a shortest $(v_k,v_{\ell})$-path
is a graph in Figure~\ref{notlt}(vi). Suppose exactly one of $v_av_k$ and 
$v_{\ell}v_c$ is an edge of $X$ and by symmetry we assume it is $v_av_k$.
Then any vertex adjacent to $v_b$ except $v_a$ is adjacent to $v_a$. So there must 
be a vertex $v_i$ adjacent to $v_a$ but not to $v_b$. Thus the subgraph of $X$ 
induced by $v_i,v_a,v_c,v_d$ and the vertices in a shortest $(v_b,v_{\ell})$-path 
is a graph in Figure~\ref{notlt}(v). Finally suppose $v_av_k$ and $v_{\ell}v_c$ are 
both edges of $X$. Then there must be a vertex $v_i$ adjacent to $v_a$ but not to 
$v_b$ and a vertex $v_j$ adjacent to $v_c$ but not to $v_d$. The subgraph induced 
by $v_i,v_a,v_b,v_c,v_d,v_j$ and the vertices in a shortest $(v_b,v_d)$-path 
is a graph in Figure~\ref{notlt}(viii). This completes the proof.
\qed

Theorem \ref{main} follows immediately from Theorems~\ref{alt}, \ref{pig}, 
\ref{cancompletetheorem} and~\ref{cannotcompletetheorem}. Because deciding whether 
a partially oriented graph is a $C_k$ ($k \geq 4$) or one of the graphs in 
Figures~\ref{pigfigure}--\ref{notlt} can be done in polynomial time, obstructions 
for acyclic local tournament orientation completions can be recognized 
in polynomial time. If a partially oriented graph cannot be completed to
an acyclic local tournament then by Proposition~\ref{obstruction} it critically
contains an obstruction and one can find an obstruction in it by deleting vertices
and replacing arcs with edges. This again can be accomplished in polynomial time.

To conclude this paper, we make a remark on partially oriented graphs which cannot 
be completed to acyclic local tournaments and are minimal with respect to only 
vertex deletions. 
Let $Y$ be such a graph, that is, $Y$ cannot be completed to an acyclic local 
tournament and for each vertex $v$ of $Y$, $Y-v$ can be completed to an acyclic 
local tournament. Since $Y$ cannot be completed to an acyclic local tournament, by 
Proposition~\ref{obstruction} it critically contains an obstruction. Since $Y$ is 
minimal with respect to vertex deletion, $Y$ is either an obstruction described in 
Theorem~\ref{main} or is obtained from an obstruction by orienting some edges.

\end{document}